\documentclass[12pt]{amsart}
\usepackage{mathptmx}
\usepackage{amsmath}
\usepackage{amssymb}
\usepackage{array}
\usepackage{geometry}
\usepackage[bookmarks=true,colorlinks=true, pdfstartview=FitV, linkcolor=black, citecolor=blue, urlcolor=black]{hyperref}
\usepackage{color}
\definecolor{DarkRed}{rgb}{0.55,.00,0.2}
\definecolor{DarkGrey}{rgb}{0.35,.35,0.35}

\theoremstyle{definition}

\theoremstyle{remark}

\numberwithin{equation}{section}

\begin{document}

\title{New inversion, convolution and Titchmarsh's  theorems for the half-Hilbert transform}

\author{S. Yakubovich}
\address{Department of Mathematics, Fac. Sciences of University of Porto,Rua do Campo Alegre,  687; 4169-007 Porto (Portugal)}
\email{ syakubov@fc.up.pt}

\keywords{Keywords  here} \subjclass[2000]{44A15, 44A35, 45E05,
45E10 }

\keywords{Hilbert transform, convolution method, Mellin transform,
Titchmarsh theorem, singular integral equations }

\maketitle

\markboth{\rm \centerline{ S.  YAKUBOVICH}}{}
\markright{\rm \centerline{THE HALF-HILBERT TRANSFORM }}

\begin{abstract}
While exploiting the generalized Parseval equality for the Mellin
transform, we derive the reciprocal inverse operator in the weighted
$L_2$-space related to the Hilbert transform on the nonnegative
half-axis. Moreover, employing the convolution method, which is
based on the Mellin-Barnes integrals, we prove the corresponding
convolution and Titchmarsh's theorems for the half-Hilbert
transform. Some applications to the solvability of a new class of
singular integral equations are demonstrated. Our technique does not
require the use of methods of the Riemann-Hilbert boundary value
problems for analytic functions. The same approach will be applied
in the forthcoming research  to invert the half-Hartley transform
and to establish its convolution theorem.
\bigskip
\end{abstract}

\section{Introduction and auxiliary results}

The main object of the present paper is the Hilbert transform on the
half-axis (the half-Hilbert transform)
$$(H_+ f)(x)\equiv h(x)=  {1\over \pi} PV \int_0^\infty {f(t)\over t-x} \ dt,\quad  x
\in \mathbb{R}_+,\eqno(1.1)$$
where $f(t)$ is a complex-valued function of the space
$L_2(\mathbb{R}_+; t^{2\alpha-1}dt), \ 0 < \alpha < 1/2$ with the
norm
$$||f||_{L_2(\mathbb{R}_+;\   x^{2\alpha-1}dx)} = \left(\int_0^\infty |f(t)|^2 t^{2\alpha-1}
dt\right)^{1/2}.\eqno(1.2)$$
Integral (1.1) is understood in the principal value sense. Mapping
and inversion properties  of the Hilbert transform with the
integration over $\mathbb{R}$ are  well-known, for instance,  in
$L_p$-spaces \cite{tit} and in connection with the solvability of
singular integral equations related to the Riemann-Hilbert boundary
value problems for analytic functions (see in \cite{gakh}).
Nevertheless the use of this theory to investigate operator (1.1)
meets certain obstacles, in particular, to find its reciprocal
inverse. The problem is indeed important in mathematical physics applications. This is  why some attempts to solve it
were undertaken, for instance,  in \cite{pip}, \cite{half}, \cite{err}, reducing to the corresponding Riemann-Hilbert boundary value problem.

Our natural approach will involve the $L_2$-theory of the Mellin
transform \cite{tit}
$$(\mathcal {M} f)(s)= f^*(s)= \int_0^\infty f(t) t^{s-1}dt, \ s \in
\sigma_\alpha =\{ s \in \mathbb{C}, s=\alpha +i\tau\},\eqno(1.3)$$
where the integral  is convergent in the mean square sense with
respect to the norm in $L_2(\sigma_\alpha)$. Reciprocally,  the
inversion formula takes place
$$f(x)= {1\over 2\pi i}\int_{\sigma_\alpha} f^*(s) x^{-s} ds,\ x >0\eqno(1.4)$$
with the convergence of the integral in the mean square sense with
respect to the norm in $L_2(\mathbb{R}_+;\\  x^{2\alpha-1}dx)$.
Furthermore, for any $f_1, f_2  \in L_2(\mathbb{R}_+; x^{2\alpha-1}dx)$ the generalized Parseval
identity holds
$$\int_0^\infty f_1\left({x\over t}\right) f_2(t) {dt\over t} = {1\over 2\pi i}\int_{\sigma_\alpha}  f_1^*(s)f_2^*(s) x^{-s}
ds, \ x >0\eqno(1.5)$$
and Parseval's equality of squares of $L_2$- norms (see (1.2))
$$\int_0^\infty |f(x)|^2 x^{2\alpha-1} dx = {1\over 2\pi}
 \int_{-\infty}^{\infty} \left|f^*\left(\alpha + i\tau\right)\right|^2 d\tau.\eqno(1.6)$$

The paper is organized as follows. In Section 2 we prove an analog of the Titchmarsh  theorem (see \cite{tit}, Th. 90) that the half-Hilbert transform (1.1) extends to a bounded invertible map $H_+:  L_2(\mathbb{R}_+;\  x^{2\alpha-1}dx) \to L_2(\mathbb{R}_+;\   x^{2\alpha-1}dx), \quad  0 < \alpha < 1/2$. In Section 3 we construct a convolution operator related to the half-Hilbert transform by a  method, based on the double Mellin-Barnes integrals and developed by the author in 1990 \cite{con}, \cite{hai}, \cite{luch}. As a conclusion the convolution theorem and Titchmarsh's theorem about the absence of divisors of zero in the convolution product  will be proved. Finally in Section 4 we consider the solvability of a new class of singular integral equations, which can be solved via the half-Hilbert transform (1.1).

As the author observed,  the same technique can be applied and will be done in the forthcoming paper to study mapping,  convolution properties and find the reciprocal inverse for the half-Hartley transform
$$g(x)= \int_{\mathbb{R}_+}  f(t)\left[\cos(xt)+ \sin(xt)\right] dt, \quad x > 0,$$
giving a rigorous motivation of these results and simplifying the inversion formula in \cite{hart}.

\section{Inversion theorem}

The following formula is well  known \cite{prud},  Vol. 1, relation (2.2.4.26)
$$PV \  {1\over \pi} \int_0^\infty {t^{s-1}\over 1-t}\ dt = \cot(\pi s), \quad   0 < {\rm Re} s < 1.\eqno(2.1)$$
Its left-hand side can be written in the form
$$PV \   {1\over \pi} \int_0^\infty {t^{s-1}\over 1-t}\ dt= \lim_{\varepsilon \to 0} \varphi_\varepsilon(s),$$
where
$$\pi \varphi_\varepsilon (s)=  \left(  \int_0^{1-\varepsilon} +  \int_{1+\varepsilon}^\infty  \right) {t^{s-1}\over 1-t}\ dt,  \quad 0 < \varepsilon < 1, \   0 < {\rm Re} s < 1.\eqno(2.3)$$

{\bf Lemma 1.} {\it Let $0< \varepsilon < 1/2,\  0 < \alpha < 1$ and $s= \alpha +i\tau, \ \tau \in \mathbb{R}.$ Then
$$ \left|\varphi_\varepsilon(s)\right| \le C |s|,\eqno(2.4)$$
where $C >0$ is an absolute constant.}

\begin{proof}   Splitting  integrals in (2.3) as follows
$$\pi \varphi_\varepsilon (s)=  \left( \int_0^{1/2}+   \int_{1/2}^{1-\varepsilon} + \int_{1+\varepsilon}^{3/2} +
 \int_{3/2}^\infty  \right) {t^{s-1}\over 1-t}\ dt$$$$=  I_1(s)+ I_2(s)+ I_3(s)+ I_4(s),$$
we observe
$$ \left|I_1(s)\right| \le  \int_0^{1/2}{t^{\alpha-1}  dt\over 1-t }= O(1).$$
 Analogously,
$$ \left|I_4(s)\right| \le  \int_{3/2}^\infty { t^{\alpha-1} dt\over t-1 }= O(1).$$
Concerning integral $I_2$, we have $(s= \alpha  +i\tau)$
$$I_2(s)=  \int_{1/2}^{1-\varepsilon}  t^{\alpha-1}  {\cos(\tau\log t) + i\sin(\tau\log t)\over 1-t}\ dt $$
and via elementary inequality $|\sin x | \le |x|, \ x \in \mathbb{R}$
$$\left| \int_{1/2}^{1-\varepsilon} t^{\alpha-1}  {\sin(\tau\log t )\over 1-t}\ dt \right| \le |\tau| \int_{1/2}^{1} t^{\alpha-1}
 {|\log t | \over 1-t}\ dt =O(\tau).$$
Further,
$$\int_{1/2}^{1-\varepsilon}  t^{\alpha -1} {\cos(\tau\log t)\over (1-t)\sqrt t}\ dt = \int_{1/2}^{1-\varepsilon} t^{\alpha-1}  {\cos(\tau\log t) -1\over 1-t}\ dt  + \int_{1/2}^{1-\varepsilon}  {t^{\alpha-1}\over 1-t}\ dt $$
and after integration by parts in the second integral, we find
$$ \int_{1/2}^{1-\varepsilon}  {t^{\alpha-1} \over 1-t}\ dt = - (1-\varepsilon)^{\alpha-1}  \log \varepsilon -   2^{1-\alpha}  \log 2 +   (\alpha-1)  \int_\varepsilon^{1/2} (1-t)^{\alpha-2} \log t \ dt.$$
In the meantime, with the Lagrange theorem
$$ {\cos(\tau\log t) -1\over t-1} =   -  \tau \   {\sin(\tau\log (\xi_t)) \over \xi_t},\quad 1/2\le   t< \xi_t< 1. $$
Hence,
$$\left| \int_{1/2}^{1-\varepsilon}  t^{\alpha-1} {\cos(\tau\log t) -1\over 1-t}\ dt \right| \le 2 |\tau| \int_{1/2}^{1}
t^{\alpha-1}  dt = O(\tau).$$
Similarly,
$$I_3(s)=  \int_{1+\varepsilon}^{3/2}  t^{\alpha-1} {\cos(\tau\log t) + i\sin(\tau\log t)\over 1-t}\ dt $$
and
$$\left| \int_{1+\varepsilon}^{3/2} t^{\alpha-1}  {\sin(\tau\log t)\over 1-t }\ dt\right|\le  |\tau| \int_{1}^{3/2} t^{\alpha -1}
 {\log t  \over t-1}\ dt =O(\tau).$$
Meanwhile,
$$\int_{1+\varepsilon}^{3/2}  t^{\alpha-1}  {\cos(\tau\log t)\over 1-t }\ dt
 = \int_{1+\varepsilon}^{3/2}  t^{\alpha-1}  {\cos(\tau\log t)-1 \over 1-t}\ dt  + \int_{1+\varepsilon}^{3/2}   {t^{\alpha-1} \over 1-t}\ dt  $$
and , in turn, with the same arguments
$$ \int_{1+\varepsilon}^{3/2} { t^{\alpha-1} \over 1-t}\ dt  =  (1+\varepsilon)^{\alpha-1} \log \varepsilon  + (2/3)^{1-\alpha}  \log 2 + (\alpha-1)  \int_{\varepsilon}^{1/2}  {\log t\over (1+t)^{2-\alpha}}\ dt,$$
$$\left| \int_{1+\varepsilon}^{3/2}  t^{\alpha-1}  {\cos(\tau\log t)-1 \over 1-t}\ dt \right| \le  |\tau| \int_{1}^{3/2}  t^{\alpha-1} dt
= O(\tau).$$
Thus,
$$\left| I_2(s)+I_3(s)\right| \le (1-\alpha) 2^{3-\alpha}  \varepsilon |\log\varepsilon| + O(1)+ O(|s|)<  2^{2 -\alpha}  \log 2 +  O(1)+ O(|s|)$$
and combining with estimates above,  we complete the proof of  inequality (2.4).
\end{proof}

Now we are ready to prove an analog of the Titchmarsh inversion theorem for the half-Hilbert transform (1.1). We have

{\bf Theorem 1}. {\it Let $0 < \alpha < 1/2.$  The half-Hilbert transform $(1.1)$ extends to a bounded invertible map
$H_+:  L_2(\mathbb{R}_+;\  x^{2\alpha-1}dx) \to L_2(\mathbb{R}_+;\   x^{2\alpha-1}dx)$ and for almost all $x \in \mathbb{R}_+$ the following reciprocal formulas
$$(H_+f)(x)\equiv h(x)= - {1\over \pi } {d\over dx} \int_0^\infty \log\left| 1- {x\over t} \right| f(t) dt, \eqno(2.5)$$
$$f(x)= {1\over \pi}  \int_0^\infty \frac{h(t)}{(\sqrt x + \sqrt t)} {dt\over \sqrt t} -  (H_+ h)(x).\eqno(2.6)$$
and  the norm inequalities  hold}
$$\left|\left| f \right|\right|_{ L_2(\mathbb{R}_+;\  x^{2\alpha-1}dx)}  \le
\left|\left|  H_+ f \right|\right|_{ L_2(\mathbb{R}_+;\  x^{2\alpha-1}dx)} \le
 \cot(\pi\alpha)  \left|\left| f \right|\right|_{ L_2(\mathbb{R}_+;\  x^{2\alpha-1}dx)},\   0 < \alpha \le {1\over 4},\eqno(2.7)$$
$$\cot(\pi\alpha)  \left|\left| f \right|\right|_{ L_2(\mathbb{R}_+;\  x^{2\alpha-1}dx)}  \le
\left|\left|  H_+ f \right|\right|_{ L_2(\mathbb{R}_+;\  x^{2\alpha-1}dx)} \le
 \left|\left| f \right|\right|_{ L_2(\mathbb{R}_+;\  x^{2\alpha-1}dx)},\   {1\over 4} \le  \alpha  < {1\over 2}.\eqno(2.8)$$

\begin{proof}  Let $f$ belong to the space $C^{(3)}_c(\mathbb{R}_+)$  of three times continuously differentiable functions
of compact support, which is dense in $ L_2(\mathbb{R}_+;\  x^{2\alpha-1}dx)$.  Then integrating by part in (1.3), we find that $s^3 f^*(s)$ is bounded on $\sigma_\alpha$ and therefore $s f^*(s) \in L_1(\sigma_\alpha),\ f^*(s) \in L_2(\sigma_\alpha)$.   Hence owing to Lemma 1 and the Lebesgue dominated convergence theorem one can pass to the limit when $\varepsilon \to 0$ under the integral sign in
$$\frac{1}{2\pi i}\int_{\sigma_\alpha}  \varphi_\varepsilon (s)  f^*(s) x^{-s} ds.$$
Consequently,  employing (2.1) and making  simple changes of variables by Fubini's theorem with  the use of (1.4), we obtain  for all $x>0$
$$\lim_{\varepsilon \to 0}  \frac{1}{2\pi i}\int_{\sigma_\alpha}  \varphi_\varepsilon (s)  f^*(s) x^{-s} ds= \frac{1}{2\pi i}\int_{\sigma_\alpha}  \cot(\pi s) f^*(s) x^{-s} ds$$
$$= \lim_{\varepsilon \to 0}  \frac{1}{2\pi^2 i}\int_{\sigma_\alpha}  \left(  \int_0^{1-\varepsilon} +  \int_{1+\varepsilon}^\infty  \right) {t^{s-1}\over 1-t}\   f^*(s) x^{-s}  dt ds $$
$$= \lim_{\varepsilon \to 0}  \frac{1}{\pi} \left(  \int_0^{1-\varepsilon} +  \int_{1+\varepsilon}^\infty  \right) {f(x/t) dt \over (1-t) t } = \lim_{\varepsilon \to 0}  \frac{1}{\pi} \left(  \int_0^{x/(1+\varepsilon)} +  \int_{x/ (1-\varepsilon) }^\infty  \right) {f(t) dt \over t-x  } $$
$$= PV \   {1\over \pi}  \int_{0 }^\infty   {f(t) dt \over t-x  } =  (H_+ f)(x).$$
Thus we proved the equality
$$(H_+ f)(x) = \frac{1}{2\pi i}\int_{\sigma_\alpha}  \cot(\pi s) f^*(s) x^{-s} ds,\eqno(2.9)$$
which is valid for any $f \in C^{(3)}_c(\mathbb{R}_+)$.  Moreover, appealing to the Parseval equality (1.6), elementary inequality $|\tanh(\pi\tau)|\le 1,\   \tau \in \mathbb{R}$ and the extreme values of the monotone function
$$\psi(u)= \frac{\cos^2(\pi\alpha) +\sin^2(\pi\alpha) u}{ \sin^2(\pi\alpha) +\cos^2(\pi\alpha) u},\quad  u \in [0,1],$$
we get the estimates
$$\left|\left| H_+ f \right|\right|_{ L_2(\mathbb{R}_+;\  x^{2\alpha-1}dx)} =\left( {1\over 2\pi}
 \int_{-\infty}^{\infty} \left|\cot\left(\pi(\alpha+i\tau) \right) f^*\left(\alpha + i\tau\right)\right|^2 d\tau\right)^{1/2} $$
$$= \left( {1\over 2\pi}  \int_{-\infty}^{\infty} \psi\left(\tanh^2(\tau)\right) \left| f^*\left(\alpha + i\tau\right)\right|^2
 d\tau\right)^{1/2} \le \left( {\cot^2(\pi\alpha)\over 2\pi}  \int_{-\infty}^{\infty} \left| f^*\left(\alpha + i\tau\right)\right|^2 d\tau\right)^{1/2} $$$$=  \cot(\pi\alpha)  \left|\left| f \right|\right|_{ L_2(\mathbb{R}_+;\  x^{2\alpha-1}dx)},\  0 < \alpha \le {1\over 4},\eqno(2.10)$$
$$\left|\left| H_+ f \right|\right|_{ L_2(\mathbb{R}_+;\  x^{2\alpha-1}dx)} = \left( {1\over 2\pi}  \int_{-\infty}^{\infty} \psi\left(\tanh^2(\tau)\right) \left| f^*\left(\alpha + i\tau\right)\right|^2  d\tau\right)^{1/2}
 \le   \left|\left| f \right|\right|_{ L_2(\mathbb{R}_+;\  x^{2\alpha-1}dx)},$$
where $  1/4 \le  \alpha <  1/2$.  On the other hand, (2.9) and (1.1) yield   $f^*(s)=\tan(\pi s) h^*(s),\  s \in \sigma_\alpha$. Therefore, analogously,
$$\left|\left| f \right|\right|_{ L_2(\mathbb{R}_+;\  x^{2\alpha-1}dx)} =  \left( {1\over 2\pi}
 \int_{-\infty}^{\infty} \left|\tan \left(\pi(\alpha+i\tau) \right) h^*\left(\alpha + i\tau\right)\right|^2 d\tau\right)^{1/2}\le   \left|\left|  H_+ f \right|\right|_{ L_2(\mathbb{R}_+;\  x^{2\alpha-1}dx)}, \  0 < \alpha \le {1\over 4},$$
$$\left|\left| f \right|\right|_{ L_2(\mathbb{R}_+;\  x^{2\alpha-1}dx)} \le  \tan(\pi\alpha) \left|\left|  H_+ f \right|\right|_{ L_2(\mathbb{R}_+;\  x^{2\alpha-1}dx)}, \   {1\over 4} \le  \alpha  < {1\over 2}.$$
This proves inequalities (2.7), (2.8).

Further, since $C^{(3)}_c(\mathbb{R}_+)$ is dense in $ L_2(\mathbb{R}_+;\  x^{2\alpha-1}dx)$, there is a unique extension of $H_+$ as an invertible continuous map $H_+:  L_2(\mathbb{R}_+;\  x^{2\alpha-1}dx) \to L_2(\mathbb{R}_+;\   x^{2\alpha-1}dx)$.
Now, let $f \in L_2(\mathbb{R}_+;\  x^{2\alpha-1}dx)$. There is a sequence $\{ f_n\},\  f_n \in  C^{(3)}_c(\mathbb{R}_+)$
such that $\left|\left| f_n- f  \right|\right|_{ L_2(\mathbb{R}_+;\  x^{2\alpha-1}dx)}  \to 0, \ n \to \infty$.  Denoting by
$$h_n(x)=  \frac{1}{2\pi  i}\int_{\sigma_\alpha}   \cot(\pi s) f_n^*(s) x^{-s} ds,\eqno(2.11)$$
we observe by virtue of (2.7), (2.8)  that  $\{h_n\}$ is a Cauchy sequence and it has a limit in  $L_2(\mathbb{R}_+;\  x^{2\alpha-1}dx)$, which we will call $h$.  Hence,  integrating through in (2.11), we have
$$\int_0^x h_n(y) dy =  \int_0^x  \left(\frac{1}{2\pi  i}\int_{\sigma_\alpha}   \cot(\pi s) f_n^*(s) y^{- s} ds\right) dy\eqno(2.12)$$
Since by the Schwarz inequality
$$\int_0^x  [ h_n(y)- h(y)] dy \le  {x^{1-\alpha}\over \sqrt{2(1-\alpha)} } \   \left|\left| h_n- h \right|\right|_{ L_2(\mathbb{R}_+;\  x^{2\alpha-1}dx)} \to 0,\ n \to \infty,$$
and in the right-hand side of (2.12) one can change the order of integration by Fubini's theorem. Then passing  to the limit when $n \to \infty$ under integral signs in (2.12) due to the Lebesgue dominated convergence theorem, we   find
$$\int_0^x h (y) dy =  \frac{1}{2\pi i }\int_{\sigma_\alpha}   \cot(\pi s) f^*(s) {x^{1- s}\over 1-s}  ds.\eqno(2.13)$$
Differentiating by $x$ in (2.13),  we come out with the equality for almost all $x>0$
$$ h (x) =  \frac{1}{2\pi i } {d\over dx} \int_{\sigma_\alpha}   \cot(\pi s) f^*(s) {x^{1- s}\over 1-s}  ds,\eqno(2.14)$$
which coincides with (2.7) for any $f \in C^{(3)}_c(\mathbb{R}_+)$.

Meanwhile,  using relation (8.4.6.8) in \cite{prud}, Vol. 3,  we have
$$ \frac{1}{2\pi i }\int_{\sigma_\alpha}   {\cot(\pi s)\over 1-s} x^{-s} ds = - {1\over \pi x} \log\left|1- x\right|. $$
Hence appealing to (1.5) in (2.14), we complete the proof of  representation (2.5).   Finally, we establish the inversion formula (2.6). To do this, we write reciprocally to (2.14) for almost all $x >0$
$$f (x) =  \frac{1}{2\pi i } {d\over dx} \int_{\sigma_\alpha}   \tan(\pi s) h^*(s) {x^{1- s}\over 1-s}  ds
=   \frac{1}{2\pi i } {d\over dx} \int_{\sigma_\alpha}   \tan^2(\pi s) \cot(\pi s) h^*(s) {x^{1- s}\over 1-s}  ds$$
$$=  \frac{1}{2\pi i } {d\over dx} \int_{\sigma_\alpha}   \left[ {1\over \cos^2(\pi s)} -1 \right]  \cot(\pi s) h^*(s) {x^{1- s}\over 1-s}  ds=  \frac{1}{\pi i } {d\over dx} \int_{\sigma_\alpha}   {h^*(s) \over \sin(2\pi s)} {x^{1- s}\over 1-s} ds $$
$$- (H_+h)(x).$$
Differentiating with respect to $x$ under the integral sign in the right-hand side of the latter equality due to the absolute and uniform convergence, we appeal to the supplement formula for Euler's gamma-functions, relation (8.4.2.5) in \cite{prud}, Vol. 3 and generalized Parseval equality to obtain
$$\frac{1}{\pi i } {d\over dx} \int_{\sigma_\alpha}   {h^*(s) \over \sin(2\pi s)} {x^{1- s}\over 1-s} ds =
{1\over \pi} \int_0^\infty \frac{h(t)}{(\sqrt x + \sqrt t)} {dt\over \sqrt t}.$$
Thus we proved (2.6) and completed the proof of Theorem 1.
\end{proof}

{\bf Corollary 1}. {\it The half-Hilbert transform $(1.1)$ is an  isometric automorphism of the space $ L_2(\mathbb{R}_+;  \\ x^{-1/2} dx)$ and reciprocal formulas $(2.5),\ (2.6)$ hold}.

\begin{proof} In fact, letting $\alpha= 1/4$ in (2.7), (2.8), we find that the map $H_+$ is isometric in $ L_2(\mathbb{R}_+;  x^{-1/2} dx)$, i.e.
$$\left|\left|  H_+ f \right|\right|_{ L_2(\mathbb{R}_+;\  x^{-1/2}dx)} =  \left|\left| f \right|\right|_{ L_2(\mathbb{R}_+;\  x^{-1/2}dx)}$$
for all $f \in L_2(\mathbb{R}_+;  x^{-1/2} dx)$.  Appealing to Theorem 1,  we complete the proof.
\end{proof}

Let us consider the half-Hilbert transform (1.1) of a complex variable $z \in D= \{ z \in \mathbb{C}: \  0< a < |z|,\   0 < \delta < |\arg z| \le \pi \}$
$$(H_+ f)(z)\equiv G(z) =  {1\over \pi}  \int_0^\infty {f(t)\over t- z} \ dt,  \quad  z \in D.\eqno(2.15)$$

{\bf Theorem 2}. {\it Let $f \in  L_2(\mathbb{R}_+;\  x^{2\alpha-1}dx),\ 0 < \alpha< 1/2.$ Then $G(z)$  is analytic in $D$. Moreover, if the Mellin transform (1.3) of $f$  satisfies the condition $s f^*(s) \in L_2(\sigma_\alpha)$, then  $G(z)$  has the following representation}
$$G(z)= \frac{1}{2\pi i}\int_{\sigma_\alpha} [ \cot(\pi s) + i] f^*(s) z^{-s} ds,\quad z \in D.\eqno(2.16)$$

\begin{proof} In fact, the integrand in (2.5) is analytic in $D$ as a function of $z$.  Hence letting $z= |z| e^{i\theta}, \ \theta= \arg z$ with the Schwarz inequality we obtain
$$|G(z)| \le {1\over \pi} ||f||_{ L_2(\mathbb{R}_+;\  x^{2\alpha-1}dx)}
\left(  \int_0^\infty {t^{1-2\alpha} \over |t- |z|e^{i\theta} |^2} \ dt\right)^{1/2}= {1\over \pi}
 ||f||_{ L_2(\mathbb{R}_+;\  x^{2\alpha-1}dx)}$$$$\times
 \left(  \int_0^\infty {t^{1-2\alpha} \over t^2 - 2t |z|\cos\theta + |z|^2} \ dt\right)^{1/2}\le {1\over \pi}
 ||f||_{ L_2(\mathbb{R}_+;\  x^{2\alpha-1}dx)}\left(  \int_0^\infty {t^{1-2\alpha} \over t^2  + a^2} \ dt\right)^{1/2}\eqno(2.17)$$
when $\pi/2 \le |\theta| \le \pi$. Meanwhile, the latter integral is
calculated in \cite{prud}, Vol. 1, relation (2.2.9.24), which gives
$$ \int_0^\infty {t^{1-2\alpha} \over t^2 - 2 t |z|\cos\theta +  |z|^2 }\ dt = |z|^{-2\alpha} \
\frac{\pi\sin((1-2\alpha)\theta)}{\sin(2\pi\alpha)\
\sin\theta},\quad 0<  |\theta| < \pi.$$
Therefore, substituting this value into the previous expression and
employing the elementary inequality $|\sin x| \le |x|,\ x \in
\mathbb{R}$,  we come out with the estimate
$$|G(z)| \le {a^{-\alpha}(1-2\alpha)^{1/2} \over \sqrt 2 [\sin(2\pi\alpha)\sin\delta]^{1/2}}
||f||_{ L_2(\mathbb{R}_+;\  x^{2\alpha-1}dx)},\  \delta < |\theta|
\le {\pi\over 2},$$
which in combination with (2.17) establishes the analyticity of $G$
in $D$, since integral (2.15) converges absolutely and uniformly in
$D$. Now, let $s f^*(s) \in L_2(\sigma_\alpha)$.  Returning to
(1.4),  we see by virtue  of the Schwarz inequality,  that the
corresponding integral converges absolutely.  Substituting this into
(2.15), we change the order of integration by Fubini's theorem
because
$$\int_0^\infty {1\over |t- z|} \int_{\sigma_\alpha}  | f^*(s) t^{-s} ds | \ dt
 = \int_0^\infty {t^{-\alpha} dt \over\left [ t^2 -  2t |z|\cos\theta + |z|^2\right]^{1/2} }
 \  \int_{\sigma_\alpha}  | f^*(s) ds | < \infty$$
for any $z \in D$.  Then taking in mind the value of the integral \cite{pip}, p.177
$$\int_0^\infty {t^{-s} dt \over t- z} = \pi z^{-s} \left[\cot(\pi s) + i \right], \quad z \in D,\eqno(2.18)$$
we arrive at (2.16), completing the proof of Theorem 2.
\end{proof}

\section{Convolution operator for the half-Hilbert transform}

In this section we will construct and study mapping properties of  the convolution, related to the transformation (1.1).   Following the general convolution method developed    for integral transforms of the Mellin convolution type  \cite{con}, \cite{hai},  \cite{luch} we have

{\bf Definition 1}. Let  $0< \alpha < 1/4$,  $f, g$ be functions from $\mathbb{R}_+$ into $\mathbb{C}$ and $f^*,\  g^*$ be their Mellin transforms $(1.3)$.  Then the function $f*g$ being defined on $\mathbb{R}_+$ by the double Mellin-Barnes integral
$$(f*g)(x)= \frac{1}{(2\pi i)^2} \int_{\sigma_\alpha}  \int_{\sigma_\alpha} \frac{\cot(\pi s)\cot(\pi w)}{\cot(\pi(s+w))}
   f^*(s)g^*(w) x^{-s-w} dsdw\eqno(3.1)$$
  is called the convolution of $f$ and $g$  related to the half-Hilbert transform $(1.1)$  (provided that it exists).

{\bf Lemma 2}. {\it Let $0< \alpha < 1/4,\  f, g$ be such that their Mellin transforms $f^*,\ g^*$ satisfy conditions $s f^*(s),\  s g^*(s)  \in  L_2(\sigma_\alpha)$.  Then convolution $(3.1)$ exists as a continuous function on $\mathbb{R}_+$,  $(f*g)(x) \in L_2(\mathbb{R}_+; x^{4\alpha-1} dx)$  and the following inequalities  hold}
$$||f*g||_{L_2(\mathbb{R}_+; x^{4\alpha-1} dx)} \le {\cot^2(\pi\alpha) \over (2\alpha)^{3/2 } \pi} \left(\int_{-\infty}^\infty   |(\alpha+i\theta) g^*(\alpha+i\theta)|^2 dt\right)^{1/2} $$$$\times  \left(\int_{-\infty}^\infty \left|(\alpha+i\tau)  f^*(\alpha+ i\tau) \right|^2 d\tau \right)^{1/2},\  0 < \alpha \le {1\over 8},\eqno(3.2)$$
$$||f*g||_{L_2(\mathbb{R}_+; x^{4\alpha-1} dx)} \le {\cot^2(\pi\alpha) \over (2\alpha)^{3/2 } \pi\cot(2\pi\alpha)} \left(\int_{-\infty}^\infty   |(\alpha+i\theta) g^*(\alpha+i\theta)|^2 dt\right)^{1/2} $$$$\times  \left(\int_{-\infty}^\infty \left|(\alpha+i\tau)  f^*(\alpha+ i\tau) \right|^2 d\tau \right)^{1/2},\  {1\over 8} \le  \alpha <  {1\over 4}.\eqno(3.3)$$
\begin{proof} In fact,  similar to the proof of Lemma 1   we obtain the estimates
$$\left| \frac{\cot(\pi s)\cot(\pi w)}{\cot(\pi(s+w))} \right| \le \cot^2(\pi \alpha),\   s,w \in \sigma_\alpha,\ 0 < \alpha \le {1\over 8},$$
$$\left| \frac{\cot(\pi s)\cot(\pi w)}{\cot(\pi(s+w))} \right| \le {\cot^2(\pi \alpha)\over \cot(2\pi\alpha)} ,\   s,w \in \sigma_\alpha,\   {1\over 8} \le  \alpha <  {1\over 4}.$$
Hence with the Schwarz inequality for double integrals and   computation of elementary integrals, we find
$$|(f*g)(x)| \le {\cot^2(\pi \alpha) \over 4\pi^2  x^{2\alpha} }  \left(\int_{\sigma_\alpha}  \left|s f^* (s) \right|^2 |ds |\right)^{1/2}
 \left(\int_{\sigma_\alpha}  \left|w  g^* (w) \right|^2 |dw |\right)^{1/2}\  \int_{-\infty}^\infty {d\tau \over \tau^2 + \alpha^2}  $$
$$= {\cot^2(\pi \alpha) \over 4\pi \alpha  x^{2\alpha} }  \left(\int_{\sigma_\alpha}  \left|s f^* (s) \right|^2 |ds |\right)^{1/2}
 \left(\int_{\sigma_\alpha}  \left|w  g^* (w) \right|^2 |dw |\right)^{1/2},\  \alpha \in \left(0,\ {1\over 8}\right],$$
$$|(f*g)(x)| \le {\cot^2(\pi \alpha) \over 4\pi\alpha   x^{2\alpha} \   \cot(2\pi\alpha)}  \left(\int_{\sigma_\alpha}  \left|s f^* (s) \right|^2 |ds |\right)^{1/2}  \left(\int_{\sigma_\alpha}  \left|w  g^* (w) \right|^2 |dw |\right)^{1/2}, \   \alpha \in \left[ {1\over 8},\ {1\over 4} \right).$$
These estimates guarantee continuity of the convolution $(f*g)(x)$ on $\mathbb{R}_+$ via the Weirstrass test of the uniform convergence of the double integral (3.1) for $x \ge x_0 >0$.    Furthermore, appealing to the
Parseval equality (1.6) and making a simple change of variables $z= s+w$ in (3.1),  we get
$$\int_0^\infty |(f*g)(x)|^2 x^{4\alpha-1} dx = {1\over 8\pi^3}  \int_{-\infty}^\infty  \left| \int_{-\infty}^\infty 
 \frac{\cot(\pi (\alpha+i(\tau-\theta))\cot(\pi (\alpha+i\theta))}{\cot(\pi(2\alpha +i\tau))}\right.$$$$\left. \times   f^*(\alpha+ i(\tau-\theta))  g^*(\alpha+i\theta)d\theta \right|^2  d\tau.$$
Hence, taking $\alpha \in (0, 1/8]$,    we employ the generalized Minkowski inequality to derive
 $$\left(\int_0^\infty |(f*g)(x)|^2 x^{4\alpha-1} dx\right)^{1/2} \le   
 {1 \over (2\pi)^{3/2 } } \int_{-\infty}^\infty   | g^*(\alpha+i\theta)| \left(\int_{-\infty}^\infty 
 \left| \frac{\cot(\pi (\alpha+i(\tau-\theta))\cot(\pi (\alpha+i\theta))}{\cot(\pi(2\alpha +i\tau))}\right.\right.$$
 $$\left.\left.\times  f^*(\alpha+ i(\tau-\theta)) \right|^2 d\tau \right)^{1/2}  d\theta \le  {\cot^2(\pi\alpha) \over (2\pi)^{3/2 } \alpha} \left(\int_{-\infty}^\infty   |(\alpha+i\theta) g^*(\alpha+i\theta)|^2 dt\right)^{1/2}  \left(\int_{-\infty}^\infty   {d\theta\over \alpha^2+\theta^2} \right)^{1/2}$$
 $$\times \left(\int_{-\infty}^\infty \left|(\alpha+i\tau)  f^*(\alpha+ i\tau) \right|^2 d\tau \right)^{1/2}=  {\cot^2(\pi\alpha) \over (2\alpha)^{3/2 } \pi} \left(\int_{-\infty}^\infty   |(\alpha+i\theta) g^*(\alpha+i\theta)|^2 dt\right)^{1/2}$$
 $$\times 
  \left(\int_{-\infty}^\infty \left|(\alpha+i\tau)  f^*(\alpha+ i\tau) \right|^2 d\tau \right)^{1/2},$$
  which proves (3.2).  Analogously we establish (3.3). 
 
  \end{proof}

 This lemma drives us to the convolution theorem for the half-Hilbert transform.  Precisely, we state

 {\bf Theorem 3}.  {\it  Let $0< \alpha < 1/4$ and $f^*,\ g^*$ be the Mellin transforms of $f, g$, respectively, satisfying  conditions $s f^*(s),\ s g^*(s)  \in  L_2(\sigma_\alpha)$.   Then the Mellin transform of the convolution $(3.1)$ \   $(\mathcal { M}  (f*g) ) (s) \in L_2(\sigma_{2\alpha})$ and is equal to
 $$(\mathcal { M}  (f*g) ) (s) =   \frac{1}{2\pi  i \  \cot(\pi s)}  \int_{\sigma_\alpha }    \cot(\pi( s-w))\cot(\pi w) f^*(s-w )g^*(w)dw,\  s \in \sigma_{2\alpha}.\eqno(3.4)$$
 Moreover,  the factorization equality  holds
$$(H_+  (f*g) )(x) =  \ (H_+ f)(x) (H_+  g)(x), \quad x >0. \eqno(3.5)$$
Besides, if $s f^*(s),\ s g^*(s)  \in  L_2(\sigma_\alpha) \cap L_1(\sigma_\alpha)$, then  for all $x >0 $ the following representation takes place}
$$(f*g)(x)=    (H_+ g)(x) f(x)  -  {1  \over \pi} \int_0^\infty f(t) {\sqrt x \  g(t) - \sqrt t\  g(x) \over t-x } \  {dt\over \sqrt t}.\eqno(3.6) $$

\begin{proof}   Formula (3.4) and condition $ (\mathcal { M}  (f*g) ) (s) \in L_2(\sigma_{2\alpha})$ follow immediately from (3.1),  Lemma 2, inversion formula for the Mellin transform (1.4)  and Parseval's  equality (1.6).  Hence employing (2.9), we  invert  the order of integration by Fubini's theorem in the obtained iterated integral.  This is indeed allowed  since, evidently, $f^*,\ g^* \in L_1(\sigma_\alpha)$ if  $s f^*(s),\ s g^*(s)  \in L_2(\sigma_\alpha)$ and the cotangent functions are bounded on $\sigma_\alpha$. So,  making  the substitution $z=s-w$, we get the chain of equalities  
$$(H_+  (f*g) )(x)= \frac{1}{(2\pi i)^2} \int_{\sigma_{2\alpha}}  x^{-s} \int_{\sigma_\alpha} \cot(\pi( s-w))\cot(\pi w)
   f^*(s-w)g^*(w)dw\ ds$$
$$= \frac{1}{(2\pi i)^2} \int_{\sigma_\alpha}  \cot(\pi w) g^*(w) \int_{\sigma_{\alpha}}  x^{-z-w} \cot(\pi z)    f^*(z) dz\ dw$$
$$= (H_+ f)(x) (H_+  g)(x), \quad x >0.$$
This  proves  (3.5).    In order to prove (3.6), we return to (3.1) and via elementary trigonometric identities write it in the form
$$(f*g)(x)= \frac{1}{(2\pi i)^2} \int_{\sigma_\alpha}  \int_{\sigma_\alpha} [\cot(\pi s)+ \cot(\pi w)]
   f^*(s)g^*(w) x^{-s-w} dsdw$$
$$+ \frac{1}{(2\pi i)^2} \int_{\sigma_\alpha}  \int_{\sigma_\alpha} \tan(\pi(s+w)) f^*(s)g^*(w) x^{-s-w} dsdw.$$
Hence under condition $s f^*(s),\ s g^*(s)  \in  L_1(\sigma_\alpha)$ with the use of (2.9), (1.4),  (2.1), (2.3),  (1.5) and the estimate  (see (2.4))  $|\varphi_\varepsilon(s+w+1/2)| \le C[ |s|+ |w| +1/2]$,  it gives  as in the proof of Theorem 1
$$(f*g)(x)=  (H_+ f) (x) g(x)  + (H_+ g)(x) f(x) $$
$$- \frac{1}{(2\pi i)^2} \int_{\sigma_\alpha}  \int_{\sigma_\alpha} \cot (\pi(s+w+1/2)) f^*(s)g^*(w) x^{-s-w} dsdw$$
$$= (H_+ f) (x) g(x)  + (H_+ g)(x) f(x) $$
$$- \lim_{\varepsilon \to 0} \frac{1}{(2\pi i)^2} \int_{\sigma_\alpha}  \int_{\sigma_\alpha} \varphi_\varepsilon(s+w+1/2) f^*(s)g^*(w) x^{-s-w} dsdw$$
$$= (H_+ f) (x) g(x)  + (H_+ g)(x) f(x)   -  {1  \over \pi}PV \int_0^\infty {f(t)g(t) \over t-x } \  \sqrt{{x\over t}}  \  dt .$$
Therefore joining  the first and the last terms in the right-hand side of the latter equality, we come out with (3.6).

\end{proof}

Finally in this section we establish an analog of the Titchmarsh theorem  about the absence of divisors of zero in the convolution  (3.1). We have

{\bf Theorem 4}. {\it Let $0< \alpha < 1/4$ and $f^*,\ g^*$ be the Mellin transforms of $f, g$, respectively, satisfying  conditions $e^{\pi  |s|}  f^*(s),\  e^{\pi  |s|}  g^*(s)  \in  L_1(\sigma_\alpha)$.   If  $(f*g)(x)= 0,  \ x >0$, then either $f(x)=0$ or $g(x) =0$ on $\mathbb{R}_+$.}

\begin{proof}  In fact, the integral 
$$F(z)=  \frac{1}{(2\pi i)^2} \int_{\sigma_\alpha}  \int_{\sigma_\alpha} \frac{\cot(\pi s)\cot(\pi w)}{\cot(\pi(s+w))}
   f^*(s)g^*(w) z^{-s-w} dsdw$$
represents an analytic function in the domain $\hat{D}= \{ z \in \mathbb{C}: |z| > a >0, \  |\arg z| < \pi\}$, since under condition of the theorem it converges uniformly in $\hat{D}$. Precisely, we have $(s= \alpha+i\tau,\  w=\alpha+i\theta,\ z^{-s-w}= |z|^{-2\alpha} e^{(\tau+\theta)\arg z }$) 
$$\int_{\sigma_\alpha}  \int_{\sigma_\alpha} \left| \frac{\cot(\pi s)\cot(\pi w)}{\cot(\pi(s+w))}
   f^*(s)g^*(w) z^{-s-w} dsdw\right|$$$$ \le C_\alpha \int_{-\infty}^\infty  \int_{-\infty}^\infty e^{\pi [  |\tau| + |\theta| ]}
  \left| f^*(\alpha+i\tau) g^*(\alpha + iw)\right|  d\tau d\theta < \infty.$$
Moreover,  (3.1) yields that $F(x)= (f*g)(x)$. Thus by virtue the uniqueness theorem for analytic functions $F(z)= (f*g)(z),
 z \in \hat{D}$. On the other hand,  employing Theorem 2 and formula (2.16),  we find the following analog of equality (3.5) for complex $z \in D$, namely
 $$ (H_+  (f*g) )(z)=  (H_+ f)(z) (H_+  g)(z)+  \frac{i }{(2\pi i)^2} \int_{\sigma_{\alpha}}   \int_{\sigma_\alpha}  \frac{\cot(\pi s)\cot(\pi w)}{\cot(\pi(s+w))}    f^*(s)g^*(w) z^{-s-w} dsdw$$
$$= (H_+ f)(z) (H_+  g)(z) + i F(z), \quad z \in D.$$
Consequently,
$$(H_+  (f*g) )(z)=  (H_+ f)(z) (H_+  g)(z)+  i (f*g)(z),\quad z \in D\cap \hat{D}.\eqno(3.7)$$
Therefore, if $(f*g)(x)=0, \ x >0,$ then via the uniqueness theorem $(f*g)(z)\equiv 0,\ z \in \hat{D}$  and (3.7) yields 
$$(H_+ f)(z) (H_+  g)(z) =0,\quad z \in D\cap \hat{D}.$$ 
Since the left-hand side of the latter equality is the product of analytic functions in  $D\cap \hat{D}$, it means that either  $(H_+ f)(z) \equiv 0$, or $(H_+ g)(z) \equiv 0$  in  $D\cap \hat{D}$. To end the proof we just appeal to the uniqueness theorem for the Stieltjes transform (cf.,  for instance, in \cite{wid}, p. 336), concluding  that either $f=0$
or $g=0$ almost everywhere on $\mathbb{R}_+$.
\end{proof}

Theorem 3 gives an idea to define the convolution (3.1) for the half-Hilbert transform in the form of equality (3.6). Our goal will be to  consider properties of the convolution in different functional classes. It concerns, for instance, the convolution $f * x^{\beta-1}$, letting $g(x)=  x^{\beta-1},\  0 < \beta  < 1$, which, in turn, does not satisfy conditions of Theorem 3.  Nevertheless, with the use of (2.1) we define this operator in the form
$$  f* x^{\beta-1}=  x^{\beta-1}\left[ f(x)\cot(\pi\beta)  -  {1 \over \pi} \int_0^\infty f(t) {(x/t)^{3/2-\beta}   -  1 \over t-x } \  dt\right].\eqno(3.8) $$

{\bf Corollary 2}. {\it Let $\beta=1/2$ and $f$ be such that $\int_0^\infty \left(f(t) /t\right) dt =0$. Then there exist divisors of zero in the convolution product $(3.8)$.}

\begin{proof}  By simple substitution in (3.8)  $\beta=1/2$  we obtain the equality $f* x^{-1/2}= 0$, which proves the result.

\end{proof}

\section{A new class of singular integral equations}

This section is devoted to an application  of convolution (3.1) to a class of integral equations,  involving the half-Hilbert transform (1.1). Namely, it concerns operator (3.8) and we will consider the solvability of the  following integral equation
$$    f(x)\cot(\pi\beta)  -  {1  \over \pi} \int_0^\infty f(t) {(x/t)^{3/2-\beta}   -  1 \over t-x } \  dt  = x^{1-\beta} h(x),\  x >0,\eqno(4.1) $$
where $h(x)$ is a given function and $f$ should be determined.  Our main result of this section will be

{\bf Theorem 5}. {\it Let $0 < \beta < 1/2,\  \beta/2 < \alpha < 1/4,\  ,\  f(x) x^{\beta-3/2}, \  x^{1-\beta} h(x) \in L_2(\mathbb{R}_+; x^{2\alpha - 1} dx)$. Then integral equation has a unique solution  $f \in L_2(\mathbb{R}_+; x^{2\alpha - 1} dx)$ written in the form
$$f(x)=  \tan(\pi\beta)  x^{1-\beta} h(x) -  \int_0^\infty K_\beta \left({x\over t} \right) h(t) t^{-\beta} dt, \eqno(4.2)$$
where
$$K_\beta (x)= \frac { 2 x^{- (\beta+ 1/2)/2}  \tan(\pi\beta) \sin  \left( { \log x  \over 2\pi} \log \left( 3 \sin(\pi\beta)  + \sqrt{ 9\sin^2(\pi\beta) - 1}\right) \right)} { \pi \sinh((\log x)  / 2\pi)  \sqrt{ 9\sin^2(\pi\beta) - 1} },\eqno(4.3)
$$
and $ \  \sin(\pi\beta) >  1/3$,
$$K_\beta (x) =  { 2 x^{- (\beta+ 1/2)/2}  \tan(\pi\beta) \sinh   \left( \gamma ( \log x)/ 2\pi \right)\over \pi \sin\gamma \  \sinh((\log x)  / 2\pi)    } ,\   3 \sin(\pi\beta) =\cos\gamma < 1,\eqno(4.4)$$
$$K_\beta (x)=  {x^{- (\beta+ 1/2)/2}  \tan(\pi\beta) \log x \over  \pi^2 \sinh((\log x)  / 2\pi) },  \   3 \sin(\pi\beta)=1.\eqno(4.5)$$
Conversely, equation (4.2) has a unique solution in $L_2(\mathbb{R}_+; x^{2\alpha - 1} dx)$ written in the form (4.1).}

\begin{proof}    Parseval's  equality (1.6) and the shift property of the Mellin transform (1.3) yield conditions 
$f^*(s),\  h^*(s-\beta+1)  \in L_2(\sigma_\alpha)$.  Hence taking the Mellin transform from both sides of equality (4.1) and employing  (2.9), where the corresponding integral converges in $L_2(\mathbb{R}_+; x^{2\alpha - 1} dx)$, we use   elementary trigonometric formulas to   obtain the following algebraic equation
$$f^*(s) \cot(\pi\beta) \frac{\sin(\pi(2s-\beta)) + 3 \sin(\pi\beta)} { \sin(\pi(2s-\beta)) +  \sin(\pi\beta)}=  h^*(s-\beta+1),\quad s \in \sigma_\alpha.$$
Clearly, the real part of the numerator in the latter fraction does not equal to zero, when $\beta/2 < \alpha < 1/4$.  Therefore,  $\sin(\pi(2s-\beta)) + 3 \sin(\pi\beta) \neq 0, \  s \in \sigma_\alpha$ and 
$$f^*(s) = \tan(\pi\beta) \frac { \sin(\pi(2s-\beta)) +  \sin(\pi\beta)} {\sin(\pi(2s-\beta)) + 3 \sin(\pi\beta)} h^*(s-\beta+1)$$
$$= \tan(\pi\beta) h^*(s-\beta+1) \left[ 1-  \frac { 2} {\sin(\pi(2s-\beta)) + 3 \sin(\pi\beta)}\right].$$
Cancelling the Mellin transform in $L_2$ due to the uniqueness property and appealing to the generalized Parseval equality  (1.5), we arrive at the solution in the form
$$f(x)=  \tan(\pi\beta)  x^{1-\beta} h(x) -  \int_0^\infty K_\beta \left({x\over t} \right) h(t) t^{-\beta} dt, $$
where
$$K_\beta (x)=  {\tan(\pi\beta) \over \pi i} \int_{\sigma_\alpha}  \frac { x^{-s} \ ds } {\sin(\pi(2s-\beta)) + 3 \sin(\pi\beta)}, \ x > 0.$$
However, the latter integral can be calculated in order to get  values of the kernel $K_\beta (x)$ for different $\beta \in (0, 1/2)$.  Indeed, substituting $s=\alpha +i\tau$ and making the change of variables $u= e^{2\pi\tau}$, we obtain 
$$K_\beta (x)= {2 x^{-\alpha} e^{i\pi(2\alpha-\beta-1/2)} \tan(\pi\beta) \over \pi^2 } \int_{0}^\infty  \frac { u^{-i (\log x)  / 2\pi } \ du } {u^2+ 6  \sin(\pi\beta) e^{i\pi(2\alpha-\beta-1/2)}  u - e^{2\pi  i (2\alpha-\beta)} }$$
$$=  {2 x^{-\alpha} e^{i\pi(2\alpha-\beta-1/2)} \tan(\pi\beta) \over \pi^2 } \lim_{\varepsilon \to 0+} 
\int_{0}^\infty  \frac { u^{-\varepsilon - i (\log x)  / 2\pi } \ du } {\left( u+  3 \sin(\pi\beta) e^{i\pi(2\alpha-\beta-1/2)}\right)^2 
+  e^{2\pi  i (2\alpha-\beta)} \left( 9\sin^2(\pi\beta) - 1\right)},$$
where the passage to the limit under the integral sign is allowed via the uniform convergence.  Hence,  in the case $\sin(\pi\beta) >  1/3$, we call (2.18) to derive
$$K_\beta (x) =  { x^{-\alpha}  \tan(\pi\beta) \over \pi^2  \sqrt{ 9\sin^2(\pi\beta) - 1} }$$$$\times  \lim_{\varepsilon \to 0+}  
\left[ \int_{0}^\infty  \frac { u^{-\varepsilon - i (\log x)  / 2\pi } \ du } {u+   e^{i\pi(2\alpha-\beta-1/2)} \left( 3 \sin(\pi\beta)
- \sqrt{ 9\sin^2(\pi\beta) - 1}\right)}\right.$$$$\left. -  \int_{0}^\infty  \frac { u^{-\varepsilon - i (\log x)  / 2\pi } \ du } {u+   e^{i\pi(2\alpha-\beta-1/2)} \left( 3 \sin(\pi\beta)
+ \sqrt{ 9\sin^2(\pi\beta) - 1}\right)}\right]$$
$$=  { x^{- (\beta+ 1/2)/2}  \tan(\pi\beta) \over \pi i \sinh((\log x)  / 2\pi)  \sqrt{ 9\sin^2(\pi\beta) - 1} } $$$$\times 
\left[ \left( 3 \sin(\pi\beta)  - \sqrt{ 9\sin^2(\pi\beta) - 1}\right)^{ - i (\log x)  / 2\pi } -  \left( 3 \sin(\pi\beta)  + \sqrt{ 9\sin^2(\pi\beta) - 1}\right)^{ - i (\log x)  / 2\pi }\right]$$
$$= { 2 x^{- (\beta+ 1/2)/2}  \tan(\pi\beta) \over \pi \sinh((\log x)  / 2\pi)  \sqrt{ 9\sin^2(\pi\beta) - 1} } 
\sin  \left( { \log x  \over 2\pi} \log \left( 3 \sin(\pi\beta)  + \sqrt{ 9\sin^2(\pi\beta) - 1}\right) \right),$$
In the same manner,  taking $\sin(\pi\beta) <   1/3$, it becomes  
$$K_\beta (x) =  { 2 x^{- (\beta+ 1/2)/2}  \tan(\pi\beta) \sinh   \left( \gamma ( \log x)/ 2\pi \right)\over \pi \sin\gamma \  \sinh((\log x)  / 2\pi)    },$$
where $\cos \gamma = 3 \sin(\pi\beta)$.  Finally, when $\sin(\pi\beta) = 1/3$, we appeal to the simple beta-integral to deduce (4.5)
$$K_\beta (x)= {2 x^{-\alpha} e^{i\pi(2\alpha-\beta-1/2)} \tan(\pi\beta) \over \pi^2 } 
\int_{0}^\infty  \frac { u^{ - i (\log x)  / 2\pi } \ du } {\left( u+   e^{i\pi(2\alpha-\beta-1/2)}\right)^2 }
= {x^{- (\beta+ 1/2)/2}  \tan(\pi\beta) \log x \over  \pi^2 \sinh((\log x)  / 2\pi) }, $$
which is just the limit case of (4.3),  (4.4) when $\gamma \to 0$.    Analogously we prove the converse statement of the theorem.

\end{proof}

\noindent {{\bf Acknowledgments}}\\
The present investigation was supported, in part,  by the "Centro de
Matem{\'a}tica" of the University of Porto.\\

\end{document}